\documentclass[a4paper,11pt]{article}

\usepackage{amssymb,amsmath,amsthm,amsfonts}
\usepackage[all]{xy}

\newtheorem{theoreme}{Theorem}[section]

\newtheorem{prop}[theoreme]{Proposition}

\newenvironment{demo}{\noindent {\sl Proof}. \ }{\qed}

\newtheorem{stheoreme}{Theorem}[subsection]
\newtheorem{sdefin}[stheoreme]{Definition}
\newtheorem{sprop}[stheoreme]{Proposition}
\newtheorem{slemme}[stheoreme]{Lemma}
\newtheorem{scorol}[stheoreme]{Corollary}
\newtheorem{sremark}[stheoreme]{Remark}

\def\pprod{\prod\limits}
\def\ssum{\sum\limits}

\def\Hoch{H\!H}

\font\tenCal=cmsy10
\newfam\Calfam
\textfont\Calfam=\tenCal

\def\qed{\hfill{$\sqcap\!\!\!\!\sqcup$}}

\def\wt{\widetilde}

\def\L{\Lambda}
\def\l{\lambda}
\def\a{\alpha}

\title {On Morita equivalence for simple Generalized Weyl algebras}

\author{Lionel Richard\footnote{University of Edinburgh, School of Mathematics and Maxwell Institute for Mathematical Sciences,
 JCMB - Mayfield Road,
  Edinburgh EH9 3JZ, United Kingdom.
 lionel.richard@ed.ac.uk} \ and Andrea Solotar\footnote{Dto. de
Matem\'atica, Facultad de Cs. Exactas y Naturales. Universidad de
Buenos Aires. Ciudad Universitaria Pab I. (1428), Buenos Aires -
Argentina. asolotar@dm.uba.ar}}

\begin{document}\maketitle
\begin{footnotetext}{
{Research partially supported by \textsc{UBACyT X169},\textsc{PIP-CONICET 5099},
\textsc{PICS-CNRS 3410},  and \textsc{Cooperaci\'on
Internacional--CONICET-CNRS}. The first author is supported by an EPSRC grant (EP/D034167/1), the second author is a research member of
\textsc{CONICET} (Argentina).} }
\end{footnotetext}

\begin{abstract}
We give a necessary condition for Morita equivalence of simple Generalized Weyl algebras
of classical type. We propose a reformulation of Hodges' result, which describes Morita equivalences in case the polynomial defining the Generalized Weyl algebra has degree $2$, in terms of isomorphisms of
quantum tori, inspired by similar considerations in noncommutative differential geometry.
We study how far this link can be generalized for $n\ge 3$.
\end{abstract}

\section{Introduction}
The aim of this paper is to describe Morita equivalence of generalized Weyl algebras of
type $k[h](\sigma_{cl},a)$, where $\sigma_{cl}(h)=h-1$ and $a\in
k[h]$ is a polynomial, under the assumption  that the
algebra  is {\sl simple} and has finite global dimension.
Generalized Weyl algebras were introduced by V.Bavula \cite{Bavrep}, and produce a common framework for the study of some classical algebras and their quantum counterpart.
 Examples of GWA are, $n$-th Weyl algebras, $\mathcal{U}(\mathfrak{sl}_2)$,
primitive quotients of $\mathcal{U}(\mathfrak{sl}_2)$, its quantized versions and also the subalgebras of invariants of
these algebras under the action of finite cyclic subgroups of automorphisms.
 It results from the discussion in \cite[\S 3.1]{JAA} that from the point of view of Morita equivalence these two cases (classical . quantum) might be treated separately. We focus here on the classical case,  also studied by T.J.Hodges \cite{H2} under the
name of {\em Non commutative deformations of type
$A$-Kleinian singularities}.

These algebras are naturally ${\mathbb Z}$-graded, and they play a crucial role in a recent paper \cite{Sierra} by Susan Sierra on rings graded equivalent to the Weyl algebra. 
Nevertheless we are dealing here with usual Morita equivalence, and the grading will not play any visible role in the following.

\medskip

{\bf Notation.} For $a\in k[h]$, denote $A(a)=k[h](\sigma_{cl},a)$
the $k$-algebra generated over $k[h]$ by two generators $x,y$
satisfying the relations
\begin{equation} \label{reldef}
\begin{array}{c}
xh=(h-1)x,\\
yh=(h+1)y,\\
yx=a(h),\\
xy=a(h-1).
\end{array}
\end{equation}

 We recall the following result, which follows from  
\cite[Proposition 2 and Corollary 2]{BaIdeals}.

\begin{prop} \label{simple}
 The classical GWA $A(a)=k[h](\sigma_{cl},a)$ is simple if
and only if  for any two distinct roots $\a$ and $\beta$ of  the
polynomial $a$, then $\a-\beta\not\in{\mathbb Z}$.
\end{prop}
\qed

\medskip

Furthermore, we will assume in the following that the polynomial
$a$ has distinct roots. Thanks to \cite[Theorem 4.4]{H2}, this is
equivalent to saying that the algebra $A(a)$ has finite global
dimension. We write explicitly this condition for further use:
\begin{equation}\label{s+fgld}
 \l_i-\l_j\not\in{\mathbb Z}, \ \ \ \forall i\neq j.
\end{equation}

 In this paper we will make use of  the proof given by
Hodges for ${\mathcal B}_{\l}$'s in \cite{H1}, which are exactly
the GWAs defined by a polynomial of degree 2, using additional
results from \cite{H2}. 
However, we propose a reformulation of Hodges' result,
relying on the link with quantum tori, inspired by similar considerations in noncommutative differential 
geometry \cite{Rieffel,Manin}.
It is natural to study, then, how far this link can be generalized for $n\ge 3$.
The paper is constructed as follows. Next section is dedicated to our main result Theorem \ref{mainthm}. Along the way we will study in deep detail the link between $K_0(A)$ and $\Hoch_0(A)$. In Section 3  we explicit our result in the case $n=3$, and investigate how far our necessary condition is to be sufficient. At last, in Section 4 we present some links with quantum tori, inspired by similar considerations in noncommutative differential geometry \cite{Rieffel,Manin}.

In all the following $k$ is an
algebraically closed field of characteristic zero, and in Sections \ref{sectdeg3} and \ref{secqtori} we will specify $k={\mathbb C}$.


\section{Framework} \label{sectframe}

\subsection{Normal form and degree 2 case}

We recall the following result of Bavula and Jordan.

\begin{stheoreme}[\cite{BJ}, Theorem 3.28]\label{thmisogwa}
 For $a_1,a_2\in k[h]$, $A(a_1)\simeq A_(a_2)$ if and only if $a_2(h)=\rho a_1(\epsilon h +\beta)$ for some $\epsilon\in\{-1,1\}$ and $\rho, \beta\in k$ with $\rho\neq 0$. 
\end{stheoreme}
\qed

Thanks to this result we will always assume our polynomials to be monic, i.e. we will write them in the form $a(h)=\prod_{i=1}^n(h-\l_i)$ with $\l_1,\ldots,\l_n$ the roots of the polynomial $a(h)$. Note that we may also, up to isomorphism, translate all roots by the same $-\beta$ and change the sign of all of them.

\begin{sremark}{\rm
 It follows from \cite{H2} that the polynomials defining two Morita equivalent GWAs must have the same degree.}
\end{sremark}

Before studying the general case we recall the following result in degree 2. 

\begin{stheoreme}[\cite{H1}, Theorem 5]\label{thmblambda}
 Let $a(h)=(h-\l_1)(h-\l_2)$ and $b(h)=(h-\mu_1)(h-\mu_2)$ be two polynomials of degree 2. Then $A(a)$ and $A(b)$ are Morita equivalent if and only if $\l_1-\l_2=\pm(\mu_1-\mu_2)+m$ for some $m\in{\mathbb Z}$.
\end{stheoreme}
\qed

\subsection{A sufficient condition}

The following is  a direct consequence of 
\cite[Lemma 2.4 and Theorem 2.3]{H2}.
\begin{prop} \label{prophodges}
Set $a,b\in k[h]$ two polynomials   with distinct
roots respectively $\{\l_i, 1\leq i\leq n\}$ and $\{\mu_i,  1\leq i\leq n\}$,
satisfying condition (\ref{s+fgld}).
Suppose that there exist $\tau\in{\mathcal S}_n$ and
$(m_1,\ldots,m_n)\in{\mathbb Z}^n$ such that $\l_i=\l'_{\tau(i)}+m_i$
for all $1\leq i\leq n$. Then the GWA's $A(a)$ and $A(b)$ are
Morita-equivalent.
\end{prop}
\qed

Note that for $n=2$ this condition is equivalent to the one appearing in Theorem \ref{thmblambda}.

\subsection{Morita equivalence and trace function}\label{subsecdiag}

 In the rest of this Section we study
necessary conditions for Morita equivalence.
Assume that $a$ and $b$ are two polynomials in $k[h]$, with
simple roots having non-integer differences, such that $A(a)$
and $A(b)$ are $k$-linearly Morita equivalent. Such an equivalence from the category of (say) left $A(a)$-modules to left
$A(b)$-modules is given by tensoring with a  
bimodule ${}_{A(b)}P_{A(a)}$, finitely generated and projective as a left and a right module. 
This functor induces a group isomorphism $K_0(F)$
between $K_0(A(a))$ and $K_0(A(b))$ and a $k$-linear isomorphism
$\Hoch_0(F)$ between $\Hoch_0(A(a))$ and $\Hoch_0(A(b))$. 
Here as usually $K_0(A)$ denotes the Grothendieck group of $A$, generated by finitely generated projective modules, and $\Hoch_0(A)$ the Hochschild homology space in degree zero, which is also the $k$-vector space of traces $A/[A,A]$. 

Moreover, $K_0(F)$
must preserve the usual rank function $rk:K_0(A)\to {\mathbb Z}$, defined on a projective $P$ as the length of the ${\rm Frac}A$-module $({\rm Frac}A)\otimes_{A}P$. So if we denote $\wt K_0(A)={\rm Ker}(rk)$,
 we have
the following commutative diagram:

\begin{equation}\label{diagram}
 \xymatrix{ {\wt K_0(A(a))} \ar[r]^i \ar[d]^{\wt K_0(F)} &
K_0(A(a)) \ar[r]^{tr} \ar[d]^{K_0(F)} & \Hoch_0(A(a))
\ar[d]^{\Hoch_0(F)}\\
{\wt K_0(A(b))} \ar[r]^i & K_0(A(b)) \ar[r]^{tr}  & \Hoch_0(A(b)) }
\end{equation}

Here $i$ denotes the canonical injection and $tr$ the usual
Hattori-Stallings trace map. Remark that $\wt K_0(F)$ is an
isomorphism of groups too (for more details see \cite{BEG}).
Following the ideas of \cite{H1}, we will describe as precisely as
possible the maps $\wt K_0(F)$ and $\Hoch_0(F)$.

\subsection{A basis for $\wt K_0(A(a))$}

Let $a\in k[h]$ be a polynomial of degree $n$ with simple roots
satisfying (\ref{s+fgld}).
Thanks to  \cite[Theorem 3.28]{BJ}, we can assume that 
$a(h)$ is monic, that is 
$a(h)=\prod_{i=1}^n(h-\l_i)$. Then
thanks to   \cite[Theorem 3.5]{H2} and Quillen's localization
sequence \cite{Hodgln1404} we know (by an argument analogous to
\cite[Proposition 1]{H1}) that $[A(a)], [P_1^{(a)}], \ldots, [P_{n-1}^{(a)}]$ form a
basis of $K_0(A(a))$, with $P_i^{(a)}=A(a)x + A(a)(h-\l_i)$. Moreover, thanks to 
\cite[Lemma 2.4]{H2}, we know that the $P_i^{(a)}$  are
progenerators, and give Morita equivalences between $A(a)$ and
$A(b_i)$, where $b_i=(h-\l_i-1)\prod_{j\neq i}(h-\l_j)$ is the polynomial obtained from $a$ by replacing 
$\l_i$ with $\l_i+1$. Then one easily verifies that

\begin{sprop}
With the notations above $([P_i^{(a)}]-[A(a)],\
1\leq i\leq n-1)$ is a basis of $\wt K_0(A(a))$.
\end{sprop}
\qed

\subsection{Trace of $[P_i^{(a)}]-[A(a)]$}\label{trace}

 We
compute here the trace of the projective $P_i^{(a)}$.

\begin{sprop}
Let $a(h)\in k[h]$ be a polynomial of degree at least 2 satisfying
the criterion of Proposition \ref{simple}, and denote $A=A(a)$.
Factorize $a(h)=u(h)w(h)$ with $u$ and $w$
non-constant polynomials. Assume that $u$ and $w$ are relatively
prime polynomials. The left $A$-ideal $P=Ax+Aw(h)$ is projective,
and its trace is the class of the polynomial
$1+w(h)B(h)-w(h-1)B(h-1)$, where $B(h), C(h)$ are two polynomials
such that $B(h)w(h)+C(h)u(h)=1$.
\end{sprop}
\begin{demo}
Consider the epimorphism of $A$-modules $G:A\oplus A\to P$ defined
by $G(1,0)=x,\ G(0,1)=w(h)$. Then one may easily check that $G$
admits a section $F: P\to A\oplus A$ defined by
$F(x)=(C(h-1)u(h-1),B(h-1)x), F(w(h))=(C(h)y,w(h)B(h))$. Then
$tr(P)$ is nothing but the usual trace of the idempotent $F\circ
G\in M_2(A)$, and one concludes using the defining relation
between $B$ and $C$.
\end{demo}

\medskip

{\bf Notations.} $\bullet$ Since $a(h)=(h-\l_1)\ldots(h-\l_n)$ has degree $n$,
 we see from 3.1.1 in \cite{FSSA} that $\Hoch_0(A(a))$ is naturally
isomorphic to the subspace of $k[h]$ spanned by the classes of $1, h, \ldots, h^{n-2}$. For convenience we will denote $1_a$ and $h^p_a$ the
classes of $1$ and $h^p$, so that $\Hoch_0(A(a))=k.1_a\oplus\bigoplus_{p=1}^{n-2} k.h^p_a$. Similarly
$\Hoch_0(A(b))=k.1_b\oplus\bigoplus_{p=1}^{n-2} k.h^p_b$.

$\bullet$ For any integer $\rho\geq 0$ denote $k_{\rho}$ the space of 
polynomials of degree not greater than $\rho$. 
Given $n$ distinct scalars $\l_1,\ldots,\l_n$, 
denote by $(v_1,\ldots,v_n)$ the basis of $k_{n-1}$ 
consisting of Lagrange interpolation polynomials associated  to 
$(\l_1,\ldots,\l_n)$, i.e. $v_i=u_i/r_i$,
with
$u_i=\prod_{j\neq i}(h-\l_j)$ and $r_i=\prod_{j\neq i}(\l_i-\l_j)=u_i(\l_i)$.

In fact, each $u_i$ is the quotient of two Vandermonde determinants,

$$\begin{array}{l}
u_i=\frac{V(\l_1, \dots, \l_{i-1},h,\l_{i+1},\dots,\l_n)}
{V(\l_1, \dots, \hat{\l_i},\dots,\l_n)}.
\Bigl(\frac{ V(\l_1, \dots,\l_n)}
{V(\l_1, \dots, \hat{\l_i},\dots,\l_n)}\Bigr)^{-1}=\
\frac{V(\l_1, \dots, \l_{i-1},h,\l_{i+1},\dots,\l_n)}
{V(\l_1, \dots,\l_n) }
\end{array}$$

with the  convention that $V(\l_1,\ldots,\l_n)$ is the determinant of the $n\times n$ matrix with $(i,j)$th entry   $\l_j^{n-i}$ for all $1\leq i,j\leq n$.

\begin{sprop} \label{formtrace}
Set $a(h)=(h-\l_1)\ldots(h-\l_n)$. Let
$P_i^{(a)}=A(a)x + A(a)(h-\l_i)$ for all $1\leq i\leq n-1$ be the left $A(a)$-modules
considered above. Then
$ tr(P_i^{(a)})$ is the class of the polynomial $1+(\sigma-1)u_i^{(a)}/r_i$.
\end{sprop}
\begin{demo}
We apply the preceding Proposition with $v(h)=a(h)$, $u(h)=\prod_{j\neq i}(h-\l_j)$ and $w(h)=h-\l_i$. The euclidian division of $u$ by $w$ gives $u=(h-\l_i)Q+r_i$, with ${\rm deg}Q=n-1$. Setting
$B(h)=-Q/r_i$ and $C(h)=1/r_i$, one gets
$tr(P_i^{(a)})$ as the class of the polynomial $1+\frac{r_i-u}{r_i}-\frac{\sigma(r_i-u)}{r_i}$. One concludes then using the fact that $\sigma$ is an algebra morphism.
\end{demo}

\medskip

Since the polynomial giving the trace of $P_i^{(a)}$ is of degree $n-2$, it may be identified with its class in $\Hoch_0(A(a))$. Denote $p_i^{(a)}$ the image of $(\sigma-1)(v_i)$ in $\Hoch_0(A(a))$, so that ${\rm tr}(P_i^{(a)})=1_a+p_i^{(a)}$.

\begin{slemme}
The set $(p_1^{(a)},\ldots,p_{n-1}^{(a)})$ is a basis of $\Hoch_0(A(a))$.
\end{slemme}
\begin{demo}
First we check that replacing $v_n(h)$ by the constant polynomial 1, the set $(v_1,\ldots,v_{n-1},1)$ is still  a basis of $k_{n-1}$.
 Let $\alpha_0, \alpha_1,\ldots,\alpha_{n-1}\in k$ such that $\alpha_0+\sum_{i=1}^{n-1}\alpha_iv_i(h)=0$. Replacing $h=\l_n$ we get $\alpha_0=0$, and we conclude thanks to the fact that the $v_i$'s are linearly independant.
 Now define the linear map  $S:k_{n-1}\to  k_{n-2}$  by $S(P)=\sigma(P)-P$. The set $(p_1^{(a)},\ldots,p_{n-1}^{(0)},0)$ is the image of the basis above by $S$.  Writing the matrix of $S$ in the canonical bases, one easily sees that it is surjective. This ends the proof.
\end{demo}

\medskip

Clearly we have the same results with $b$ instead of $a$ and the $\mu_i$'s instead of the $\l_i$'s.
We give now an interpretation of the trace polynomials $p_i^{(a)}$ in terms of Schur polynomials. 


\begin{sprop} \label{formtrace3}
Set  as before $a(h)=\prod_{i=1}^n(h-\l_i)$. 
Let
$P_i^{(a)}=A(a)x + A(a)(h-\l_i)$ for $1\le i \le n-1$.
Then
$$ tr(P_i^{(a)})=1_a+p_i^{(a)}$$
with
$$\begin{array}{l}
p_i^{(a)}=\ssum_{i=1}^n(-1)^{i+l}((h-1)^{n-l}- h^{n-l})
\sigma_{\small{({\underbrace{1,\dots, 1}_{l-1},\underbrace{0,\dots,0}_{n-l}})}}(\l_1,\dots,\hat{\l_i},\dots, \l_n).
\frac{V(\l_1, \dots, \hat{\l_i},\dots,\l_n)}{V(\l_1, \dots,\l_n)}
\end{array}$$
where
$\sigma_{\small{({\underbrace{1,\dots, 1}_{l-1},\underbrace{0,\dots,0}_{n-l}})}}(\l_1,\dots,\hat{\l_i},\dots, \l_n)$
denotes the Schur polynomial associated to the partition 
$(\underbrace{1,\dots, 1}_{l-1},\underbrace{0,\dots,0}_{n-l})$
evaluated in $(\l_1,\dots,\hat{\l_i},\dots, \l_n)$.
\end{sprop}
\begin{demo}
Recall from Proposition \ref{formtrace} that


$$p_i^{(a)}=(u_i(h-1)-u_i(h)).\frac{1}{r_i}=$$

$$
\bigl( \pprod_{j\neq i}(h-1- \l_j) -  \pprod_{j\neq i}(h- \l_j)\bigr).\frac{1}{ \pprod_{j\neq i}(\l_i- \l_j)}=$$
$$
\bigl(\frac{V(\l_1, \dots, \l_{i-1},h-1,\l_{i+1},\dots,\l_n)}
{V(\l_1, \dots, \hat{\l_i},\dots,\l_n)}- \frac{V(\l_1, \dots, \l_{i-1},h,\l_{i+1},\dots,\l_n)}
{V(\l_1, \dots, \hat{\l_i},\dots,\l_n)} \bigr).
\frac{V(\l_1, \dots, \hat{\l_i}\dots,\l_n)}
{V(\l_1, \dots,\l_n)}=$$
$$
\frac{V(\l_1, \dots, \l_{i-1},h-1,\l_{i+1},\dots,\l_n)- V(\l_1, \dots, \l_{i-1},h,\l_{i+1},\dots,\l_n)}
{V(\l_1, \dots,\l_n)}=$$
$$
\frac{\det \left( 
\begin{array}{cccccc}
	\l_1^{n-1} & \l_2^{n-1} & \cdots & (h-1)^{n-1}-h^{n-1}& \cdots &\l_n^{n-1}\\
	\vdots & \vdots & \vdots & \ddots & \vdots & \vdots \\
	\l_1 & \l_2 & \cdots & h-1-h& \cdots &\l_n \\
	1 & 1  & \cdots & \underbrace{1-1}_i &1 & 1 \\
\end{array}
 \right)}{V(\l_1, \dots,\l_n)}.
$$

\bigskip

Developing by the $i$-th column we obtain:

\bigskip

$\ssum_{i=1}^n(-1)^{i+l}((h-1)^{n-l}- h^{n-l}).
\frac{\det \left( 
\begin{array}{ccccccc}
	\l_1^{n-1} & \l_2^{n-1} & \cdots & \l_{i-1}^{n-1}&\l_{i+1}^{n-1}& \cdots &\l_n^{n-1}\\
	\vdots & \vdots & \vdots & \ddots & \vdots & \vdots \\
	\l_1 & \l_2 & \cdots & \l_{i-1}& \l_{i+1}&\cdots &\l_n \\
	1 & 1  & \cdots & 1 &1&\cdots & 1 \\
\end{array}
 \right)}{V(\l_1, \dots,\l_n)}=$
 
$\bigl(\ssum_{i=1}^n(-1)^{i+l}((h-1)^{n-l}- h^{n-l}).\sigma_{\small{({\underbrace{1,\dots, 1}_{l-1},\underbrace{0,\dots,0}_{n-l}})}}(\l_1,\dots,\hat{\l_i},\dots, \l_n). \frac{V(\l_1, \dots, \hat{\l_i}\dots,\l_n)}
{V(\l_1, \dots,\l_n)} \bigr).$
\bigskip
\end{demo}

\medskip

Let us remark here that Schur polynomials also play a central role in the classification up to Morita equivalence of Cherednik algebras in \cite{BEG}.

\subsection{Computing $\Hoch_0(F)$}

In this subsection we consider two polynomials
$a(h)=\prod_{i=1}^n(h-\l_i)$ and
$b(h)=\prod_{j=1}^n(h-\mu_j)$ with all distinct roots with non-integer differences. Assuming that the algebras $A(a)$ and $A(b)$ are Morita equivalent, and using the notations of \ref{subsecdiag},
we describe
now
$\Hoch_0(F)$ as a matrix $(\a_{ij})\in GL_n(k)$, in the bases
$(p_1^{(a)},\ldots,p_{n-1}^{(a)})$, $(p_1^{(b)},\ldots,p_{n-1}^{(b)})$.

\subsubsection{Notations} \label{notat}

$\bullet$ Set $P={}_{A(b)}P_{A(a)}$  the progenerator such that
$F\equiv P\otimes_{A(a)}(\ )$. It must have rank 1 as an $A(b)$-module
(because both rings are noetherian domains), so $[P]-[A(b)]$ has
rank 0 in $K_0(A(b))$, and there exist $m_1,\ldots,m_{n-1}\in{\mathbb Z}$ such that  $[P]=[A(b)]+m_1\big([P_1^{(b)}]-[A(b)]\big)+\ldots+m_{n-1}\big([P_{n-1}^{(b)}]-[A(b)]\big)$  in $K_0(A(b))$  and $tr^{(b)}(P)=1_b +m_1
p_1^{(b)}+\ldots+m_{n-1}p_{n-1}^{(b)}$ in $\Hoch_0(A(b))$.

$\bullet$ Because $\wt K_0(F)$ is a group isomorphism, there
exists a matrix $N=( n_{ij})\in GL_{n-1}({\mathbb Z})$ such that for all $1\leq i\leq n-1$ we have
$$\wt K_0(F)\big([P_i^{(a)}]-[A(a)]\big)=n_{1,i}\big([P_1^{(b)}]-[A(b)]\big)+\ldots+n_{n-1,i}\big([P_{n-1}^{(b)}]-[A(b)]\big).$$

It results from the definition of the $m_i$'s that the matrix associated to $K_0(F)$ with respect to the bases $\big([P_1^{(a)}]-[A(a)],\ldots,[P_{n-1}^{(a)}]-[A(a)],[A(a)]\big)$ and $\big([P_1^{(b)}]-[A(b)],\ldots,[P_{n-1}^{(b)}]-[A(b)],[A(b)]\big)$  is 
$$\left(\begin{array}{cc} N&{\begin{array}{c}m_1\\ \vdots\\ m_{n-1}\end{array}}\\ 0&1\end{array}\right).$$

\subsubsection{Link between the matrices of $\wt K_0(F)$ and $\Hoch_0(F)$} \label{parlink}

We still consider the commutative diagram (\ref{diagram}). Then
we get for all $1\leq i\leq n-1$
$$\Hoch_0(F)\big({\rm tr}([P_i^{(a)}]-[A(a)])\big)={\rm tr}\big(\wt
K_0(F)([P_i^{(a)}]-[A(a)])\big),$$
 that is
$\Hoch_0(F)(p_i^{(a)})={\rm tr}\big(n_{1,i}([P_1^{(b)}]-[A(b)])+\ldots+n_{n-1,i}([P_{n-1}^{(b)}]-[A(b)])\big)$,
so $\alpha_{1,i} p_1^{(b)}+\ldots+\alpha_{n-1,i} p_{n-1}^{(b)}=n_{1,i} p_1^{(b)}+\ldots+
n_{n-1,i}p_{n-1}^{(b)}$. Since the $p_i^{(b)}$'s are linearly independant, we get $\alpha_{k,i}=n_{k,i}$ for all $1\leq k ,i\leq n-1$, that is, the matrices associated to $\wt K_0(F)$ and $\Hoch_0(F)$ in our chosen bases are equal.

\subsubsection{Computing $\Hoch_0(F)$}

 Because the diagram (\ref{diagram}) is commutative,
we have for all $1\leq i\leq n-1$:
\begin{equation} \label{comm}
\Hoch_0(F)\big({\rm tr}([P_i^{(a)}])\big)={\rm tr}\big(K_0(F)([P_i^{(a)}])\big).
\end{equation}
 The left part of this equation is equal to
$\Hoch_0(F)(1_a+p_i^{(a)})$.

\begin{slemme}
The following formulas hold respectively in $\Hoch_0(A(a))$ and $\Hoch_0(A(b))$
\begin{equation}
1_a=-\sum_{i=1}^{n-1}(\l_i-\l_n)p_i^{(a)};\ \ \
1_b=-\sum_{j=1}^{n-1}(\mu_j-\mu_n)p_j^{(b)}. 
\end{equation} 
\end{slemme}
\begin{demo}
 We give the proof for $a(h)$, the proof for $b(h)$ being completely similar. So we omit the upper indices $(a)$ in the following. Recall from the notations introduced in \ref{trace} that $p_i(h)=(\sigma-1)(v_i(h))$, with $v_1(h),\ldots,v_n(h)$ the Lagrange interpolation polynomials associated to $\l_1,\ldots,\l_n$. Reasoning in $k_{n-1}$, we have $h=\sum_{i=1}^{n}\l_iv_i(h)$ and $1=\sum_{i=1}^nv_i(h)$, so that $h=\sum_{i=1}^{n-1}(\l_i-\l_n)v_i(h)+\l_n$. We conclude by noticing that $1=-(\sigma-1)(h)$.
\end{demo}

\medskip

Now we have
$$\begin{array}{c}
\Hoch_0(F)(1_a+p_i^{(a)})=\Hoch_0(F)\left(\ssum_{j=1}^{n-1}(-\l_j+\l_n+\delta_{ij})p_j^{(a)}\right)\\
=\ssum_{j=1}^{n-1}\ssum_{k=1}^{n-1}\big((-\l_j+\l_n+\delta_{ij})\a_{kj}p_k^{(b)}\big)
.\end{array}$$

On the other hand, we have
$$\begin{array}{c}
{\rm tr}\big(K_0(F)([P_i^{(a)}])\big)={\rm tr}\big([P\otimes_{A(a)}
P_i^{(a)}]\big)=\\
{\rm tr}\left(\ssum_{k=1}^{n-1}n_{k,i}([P_k^{(b)}]-[A(b)])+[{}_{A(b)}P]\right)
=\\
\ssum_{k=1}^{n-1}n_{k,i}p_k^{(b)}+1_b +\ssum_{k=1}^{n-1}m_k p_k^{(b)}=\\
\ssum_{k=1}^{n-1}(n_{k,i}+m_k+(-\mu_k+\mu_n))p_k^{(b)}.\end{array}$$

So Equation  (\ref{comm}) gives rise for all $1\leq k\leq n-1$ to
$$\sum_{j=1}^{n-1}(-\l_j+\l_n+\delta_{ij})\a_{kj}=n_{ki}+m_k+(-\mu_k+\mu_n).
$$

Thanks to \S \ref{parlink} we can rewrite the preceding equation only in terms of the
$n_i$'s, and finally summarize the results of this section in the
following

\begin{stheoreme} \label{mainthm}
Set $a=(h-\l_1)\ldots(h-\l_n), b=(h-\mu_1)\ldots(h-\mu_n)\in
k[h]$ two polynomials of degree $n$ such that $\l_i-\l_j\not\in{\mathbb Z}$, $\mu_i-\mu_j\not\in{\mathbb Z}$ for all $i\neq j$. Define the following column vectors: $\L=(\l_n-\l_1,\ldots,\l_n-\l_{n-1})^t, \Omega=(\mu_n-\mu_1,\ldots,\mu_n-\mu_{n-1})^t\in k^{n-1}$.
Assume the algebras $A(a)$ and
$A(b)$ are Morita equivalent. Then there exist a matrix $N=(n_{ij})\in GL_{n-1}({\mathbb Z})$ and a column vector of integers $M=(m_1,\ldots,m_{n-1})^t\in{\mathbb Z}^{n-1}$ such that:
\begin{equation}\label{cn}
N.\L=\Omega+M
\end{equation}
\end{stheoreme}
\begin{demo}
It results from the preceding computations that  for all $1\leq i,k\leq n-1$ we have the following equation
$$\sum_{j=1}^{n-1}(-\l_j+\l_n+\delta_{ij})n_{kj}=n_{ki}+m_k+(-\mu_k+\mu_n).
$$
The term $n_{ki}$ appears once on both sides of this equaliy, so cancels, and $i$ does not appear anymore in the equation. Then the statement of the theorem is just a rephrasing of these facts in terms of matrices.
\end{demo}

\begin{sremark}{\rm \begin{itemize}
 \item Since $GL_1({\mathbb Z})=\{1,-1\}$,  condition (\ref{cn}) can be considered as an extension in degree $n$ of the condition obtained by Hodges in \cite{H1} (see Theorem \ref{thmblambda}).
 \item
Condition (\ref{cn})  is actually saying that the ${\mathbb Z}$-lattice generated in $k$ by the $\l_n-\l_i$'s has to be the same as the one generated by the $\mu_n-\mu_j$'s. There is a canonical way to associate a noncommutative torus to a lattice (see \cite{Rieffel,Manin}), and we will discuss this in Section \ref{secqtori}. 
\end{itemize}}
\end{sremark}

\section{Discussion on the case of degree 3} \label{sectdeg3}

In this section and the following one we assume that $k={\mathbb C}$.
Consider two polynomials
$a(h)=(h-\l_1)(h-\l_2)(h-\l_3)$ and
$b(h)=(h-\mu_1)(h-\mu_2)(h-\mu_3)$ both
satisfying the criterion (\ref{s+fgld}). 

\subsection{Notations} \label{notat3}

$\bullet$ Set $P={}_{A(b)}P_{A(a)}$ as in the previous section. We already know that
 $[P]=[A(b)]+m_1([P_1^{(b)}]-[A(b)])+m_2([P_2^{(b)}]-[A(b)])$ and ${\rm tr}^{(b)}(P)=1_b +m_1
p_1^{(b)}+m_2p_2^{(b)}$ for some $m_1,m_2\in{\mathbb Z}$, and that there exists a matrix $N=\left(\begin{array}{cc} n_1 &n_2\\ n_3&n_4
\end{array}\right)\in GL_2({\mathbb Z})$ such that
$$\begin{array}{c}
\wt K_0(F)([P_1^{(a)}]-[A(a)])=n_1([P_1^{(b)}]-[A(b)])+n_3([P_2^{(b)}]-[A(b)])\\
\wt K_0(F)([P_2^{(a)}]-[A(a)])=n_2([P_1^{(b)}]-[A(b)])+n_4([P_2^{(b)}]-[A(b)]).
\end{array}$$

Theorem \ref{mainthm} translates in the following way in the present setting.

\begin{sprop} \label{prop3}
Set $a=(h-\l_1)(h-\l_2)(h-\l_3), b=(h-\mu_1)(h-\mu_2)(h-\mu_3)\in
k[h]$ two polynomials of degree 3  such that $\l_i-\l_j\not\in{\mathbb Z}$, $\mu_i-\mu_j\not\in{\mathbb Z}$ for all $i\neq j$. Assume the algebras $A(a)$ and
$A(b)$ are Morita equivalent. Then there exist a matrix $\left(\begin{array}{cc} n_1 &n_2\\ n_3&n_4
\end{array}\right)\in M_2({\mathbb Z})$ and integers $m_1,m_2\in{\mathbb Z}$ such that
\begin{equation}\label{cn3}
\begin{array}{l}
n_1n_4-n_2n_3=\pm 1\\
(-\l_1+\l_3)n_1+(-\l_2+\l_3)n_2=m_1+(-\mu_1+\mu_3)\\
(-\l_1+\l_3)n_3+(-\l_2+\l_3)n_4=m_2+(-\mu_2+\mu_3)
\end{array}
\end{equation}
\end{sprop}
\qed

\medskip

We shall  note that in the ``generic'' case, knowing $\l_i$'s,  $\mu_j$'s and $m_k$'s satisfying (\ref{cn}), the matrix $N$ is uniquely determined. 
More precisely, given $\l_1,\l_2,\l_3$, $\mu_1, \mu_2,\mu_3$, $m_1,m_2$ and 2 matrices $N$ and $N'$ satisfying (\ref{cn}), assume that 

\begin{equation}\label{qsimple2}
(\l_3-\l_1)/(\l_3-\l_2)\not\in{\mathbb Q}.
\end{equation}
 
 Since the vector $((\l_3-\l_1)/(\l_3-\l_2),1)$ should be in the kernel of the matrix $N-N'$, this matrix has to be null, that is $N=N'$.

\subsection{Reduction of the matrix $\Hoch_0(F)$.}\label{examplematrix}

We present in this section the matrices $\Hoch_0(F)$ associated to some elementary operations on the roots of the polynomial $a(h)$.

\subsubsection{Exchanging $\l_1$ and $\l_2$.}

We consider the polynomial $b(h)=(h-\l_2)(h-\l_1)(h-\l_3)$, that
is we set $\mu_1=\l_2,\mu_2=\l_1$ and $\mu_3=\l_3$. Obviously
$A(a)=A(b)=A$, and the Morita equivalence may be given by
$P={}_{A(b)}A(b)_{A(a)}$. Then ${\rm Tr}^{(b)}(P)=1$, and
$m_1=m_2=0$. Also we have $K_0(F)={\rm Id}$, and
$P_1^{(a)}=Ax+A(h-\l_1)=Ax+A(h-\mu_2)=P_2^{(b)}$, so that
$n_1=0,n_3=1$. Then equations (\ref{cn3}) lead to $n_2=1,n_4=0$,
and we finally
get $\left(\begin{array}{cc} n_1&n_2\\
n_3&n_4\end{array}\right)=\left(\begin{array}{cc} 0&1\\
1&0\end{array}\right)=N_1$.

\subsubsection{Exchanging $\l_2$ and $\l_3$.}

We consider the polynomial $b(h)=(h-\l_1)(h-\l_3)(h-\l_2)$, that
is we set $\mu_1=\l_1,\mu_2=\l_3$ and $\mu_3=\l_2$. Once again
$A(a)=A(b)=A$, and the Morita equivalence may be given by
$P={}_{A(b)}A(b)_{A(a)}$, so that ${\rm Tr}^{(b)}(P)=1$, and
$m_1=m_2=0$. We have  $K_0(F)={\rm Id}$, and
$P_1^{(a)}=Ax+A(h-\l_1)=P_1^{(b)}$, so $n_1=1,n_3=0$. Then
equations (\ref{cn3}) lead to $n_2=-1,n_4=-1$, and we finally
get $\left(\begin{array}{cc} n_1&n_2\\
n_3&n_4\end{array}\right)=\left(\begin{array}{cc} 1&-1\\
0&-1\end{array}\right)=N_2$.

\subsubsection{$\l_1\mapsto\l_1+1$}

By  \cite[Theorem 2.3 and Lemma 2.4]{H2}, $P=P_1^{(b)}$
provides a Morita equivalence between $A(a)$ and $A(b)$, with
$\mu_1=\l_1+1,\mu_2=\l_2,\mu_3=\l_3$. By definition of $P$, we get
$m_1=1,m_2=0$. Then the identity matrix $I_2$ satisfies the equations (\ref{cn3}).

\subsubsection{$\l_i\mapsto -\l_i+1$}

It is known after \cite{BJ} that $A(a)$ is isomorphic to $A(b)$
for $b(h)=a(1-h)$. So once again using $P={}_{A(b)}A(b)_{A(a)}$ in
this context we get $m_1=m_2=0$. The matrix $-I_2$ satisfies the equations (\ref{cn3}).
Moreover the isomorphism is given
by $x\mapsto y,y\mapsto x,h\mapsto 1-h$. 


\subsection{A subgroup of $SL_2({\mathbb Z})$}\label{subsec:c6}

The necessary condition appearing in Proposition \ref{prop3} is still weaker than the sufficient condition of Proposition \ref{prophodges}. In the following we show that the necessary condition (\ref{cn3}) cannot be sufficient in degree 3, at least not without the extra assumption that the polynomials $a$ and $b$ both satisfy (\ref{s+fgld}).

Given two polynomials $a$ and $b$, a permutation of the first two roots of
$b$ leads to  multiplication on the right of $\Hoch_0(F)$
by the matrix $N_1$. Thanks to this, we may assume
that $\Hoch_0(F)\in SL_2({\mathbb Z})$, that is $n_1n_4-n_2n_3=1$.

\medskip

{\bf Notation.} Let $G$ be  the subgroup consisting of matrices $N\in SL_2({\mathbb Z})$, such that for \emph{all} triples $(\l_1,\l_2,\l_3)$ and $(\mu_1,\mu_2,\mu_3)$ satisfying (\ref{cn3}), the algebras $A(a)$ and $A(b)$ are Morita equivalent, with $a=(h-\l_1)(h-\l_2)(h-\l_3)$ and $b=(h-\mu_1)(h-\mu_2)(h-\mu_3)$.

It is clear from paragraph \ref{examplematrix} that $-I_2$
and $N_1N_2$ belong to $G$. These two elements generate a subgroup $G_6$ isomorphic to ${\mathbb Z}/2{\mathbb Z}\times{\mathbb Z}/3{\mathbb Z}$. 
The 6 elements of this subgroup are the identity matrix $I_2$, its opposite $-I_2=\left(\begin{array}{cc}
-1&0\\
0&-1
\end{array}\right)$, $N_1N_2=\left(\begin{array}{cc}
0&-1\\
1&-1
\end{array}\right)$, $-N_1N_2=\left(\begin{array}{cc}
0&1\\
-1&1
\end{array}\right)$, $(N_1N_2)^2=\left(\begin{array}{cc}
-1&1\\
-1&0
\end{array}\right)=N_2N_1$ and $-N_2N_1=\left(\begin{array}{cc}
1&-1\\
1&0
\end{array}\right)$.

\begin{sprop}
The matrices $N_1N_2$ and $-I_2$ generate $G$; that is:
$G=G_6$.
\end{sprop}
\begin{demo}
 Let $N=\left(\begin{array}{cc}
 n_1&n_2\\
 n_3&n_4
\end{array}\right)
$ be an element of $SL_2({\mathbb Z})$. We will show that if $N$ is not one of the 6 matrices above, then there exist triples $(\l_1,\l_2,\l_3)$ and $(\mu_1,\mu_2,\mu_3)$ satisfying (\ref{cn3}), such that the algebra $A(a)$ is simple with finite global dimension and the algebra $A(b)$ is not, with $a=(h-\l_1)(h-\l_2)(h-\l_3)$ and $b=(h-\mu_1)(h-\mu_2)(h-\mu_3)$. So $N\not\in G$.

$\bullet$ Assume $|n_1n_2|>1$. Since $N\in SL_2({\mathbb Z})$, this implies $n_1\neq n_2$. Consider now the triple $\l_1=1/(2n_1)$, $\l_2=1/(2n_2)$ and $\l_3=0$. It results from the hypothesis that $0<|\l_i-\l_j|<1$ for all $i\neq j$. So the algebra $A(a)$ is simple and of finite global dimension. But for a triple $\mu_1,\mu_2,\mu_3$ satisfying (\ref{cn3}) we get $\mu_3-\mu_1=-m_1-1\in{\mathbb Z}$, so the algebra $A(b)$ is not simple, or not of finite global dimension if $m_1=-1$.

$\bullet$ The case $|n_3n_4|>1$ is dealt with similarly.

So a matrix in the group $G$ has all its entries in the set $\{0,1,-1\}$.

$\bullet$ $n_1=0$. Then necessarily $n_2n_3=-1$. Assume first that $n_2=-n_3=1$. If $n_4=0$ then denote $x=\left(\begin{array}{cc}
 0&1\\
 -1&0
\end{array}\right)$ the corresponding matrix.
Considering the triples $(\l_1=3/4+i, \l_2=1/4-i, \l_3=0)$ and $(\mu_1=-1/4+i, \mu_2=3/4+i,\mu_3=0)$ leads as before to a simple and a non simple algebras. If $n_4=-1$ then consider $(\l_1=3/4, \l_2=1/4, \l_3=0)$ and $(\mu_1=1/4, \mu_2=-1,\mu_3=0)$. Note that a similar example will do as soon as $n_3n_4=1$, or by symmetry of the problem as soon as $n_1n_2=1$, and that none of the matrices in $G_6$ satisfies such an hypothesis. At last, taking $n_4=1$ gives $N=-N_1N_2$, which belongs to $G_6$.
If $n_2=-n_3=-1$ then multiplying by $-I_2$ leads to similar conclusions. So $x\not\in G$.

$\bullet$ $n_1=1$. We consider three subcases, depending on the value of $n_2$. \begin{enumerate}
 \item $n_2=0$. Then necessarily $n_4=1$. So if $n_3=0$ then $N={\rm Id}$; if $n_3=1$ then we are in the case $n_3n_4=1$ which can be dealt with as before; if $n_3=-1$ then one can easily check $N=\left(\begin{array}{cc}
 1&0\\
 -1&1
\end{array}
\right)=x^{-1}N_2N_1$. So $N\not\in G_6$, otherwise we would have $x\in G$.
\item $n_2=1$. Then $n_1n_2=1$, and we already noticed that none of the matrices satisfying such an  hypothesis is in $G$.
\item $n_2=-1$. Then $n_4+n_3=1$, i.e. $(n_3,n_4)\in\{(1,0),(0,1)\}$. The first case corresponds to $-N_2N_1$, which belongs to $G_6$. One checks easily that the second case corresponds to the matrix $N=x(N_2N_1)^{-1}$, which cannot belong to $G$, otherwise we would have $x\in G$.
\end{enumerate}

$\bullet$ $n_1=-1$. This case is strictly similar to the preceding one, up to multiplication by the matrix $-I_2$ which belongs to $G_6$.
\end{demo}

\section{Links with  quantum tori}\label{secqtori}

As for the previous section, we assume here that $k={\mathbb C}$.

\subsection{Quantum tori}

We recall the following
\begin{sdefin}
Let $n\geq 1$ be an integer and $Q=(q_{ij})\in M_n({\mathbb C}^*)$ be a multiplicatively antisymmetric matrix (i.e. $q_{ij}q_{ji}=q_{ii}=1 \ \forall i,j$). The quantum torus (or MacConnell-Pettit algebra \cite{MCP}) parametrized by $Q$ is the ${\mathbb C}$-algebra generated by $X_1,\ldots,X_n$, with relations $X_iX_j=q_{ij}X_jX_i$, and their inverses $X_1^{-1},\ldots,X_n^{-1}$. It is denoted $T_Q={\mathbb C}_Q[X_1^{\pm 1},\ldots,X_n^{\pm 1}]$.
\end{sdefin}

These algebras play a crucial role in quantum algebra (see for example \cite{Cauchon}), and have been extensively studied. 
Note that when $n=2$ the matrix $Q$ is uniquely determined by the entry $q=q_{12}$. In this case we may denote the associated quantum torus by $T_q$ or ${\mathbb C}_q[X_1^{\pm 1},X_2^{\pm 2}]$.
We will focus in the sequel on the following property.

\begin{sprop}[\cite{MCP}, Proposition 1.3]\label{thmqsimple}
Let $T_Q={\mathbb C}_Q[X_1^{\pm 1},\ldots,X_n^{\pm 1}]$ be a quantum torus. The following conditions are equivalent:
\begin{enumerate}
\item the centre of $T_Q$ is reduced to ${\mathbb C}$;
\item $T_Q$ is a simple ring;
\item if $(m_1,\ldots,m_n)\in{\mathbb Z}^n$ satisfies
\begin{equation}\label{eqqsimple}
\prod_{k=1}^nq_{kj}^{m_k}=1 ,\ \ \forall 1\leq j\leq n
\end{equation}
then $m_i=0$ for all $1\leq i\leq n$.
\end{enumerate}
\end{sprop}
\qed

\medskip

If $n=2$ then this condition is equivalent to saying that $q$ is not a root of unity.
Since we are dealing with Morita equivalence, we may mention also the following consequence of \cite[Theorem 1.4]{MCP}, \cite[Th\'eor\`eme 4.2]{endotq} and \cite[Lemma 3.1.1]{JAA}.

\begin{stheoreme} \label{thmqiso}
 Let $Q=(q_{ij}), Q'=(q'_{ij}) \in M_n({\mathbb C}^*)$ be multiplicatively antisymmetric matrices. Assume that the quantum tori $T_Q$ and $T_{Q'}$ parametrized by $Q$ and $Q'$ are simple. Then the following are equivalent
 \begin{enumerate}
 \item $T_Q$ and $T_{Q'}$ are isomorphic;
 \item there exists $M=(m_{ij})\in GL_n({\mathbb Z})$ such that for all $i,j$ one has 
 $$q'_{ij}=\prod_{t,k}q_{kt}^{m_{ki}m_{tj}};$$
 \item $T_Q$ and $T_{Q'}$ are birationnally equivalent (i.e. have isomorphic skew-fields of fractions);
 \item $T_Q$ and $T_{Q'}$ are Morita equivalent.
\end{enumerate}
\end{stheoreme}
\qed

\medskip

If $n=2$ then condition 2. is easily seen to be equivalent to $q'=q$ or $q^{-1}$.
Now we will explain how this is related to GWAs. The next subsection is devoted to the case $n=2$.

\subsection{Rank 2 case}

Our motivation here is the survey paper \cite{Manin} by Yuri Manin. Even if  the author is there interested in differential non commutative geometry and considers smooth and rapidly decreasing functions, we will keep an algebraic point of view and only consider noncommutative Laurent polynomials.

Consider a lattice of rank two ${\mathbb Z}\oplus \theta{\mathbb Z} \subset {\mathbb C}$, with $\theta\in{\mathbb C}\setminus{\mathbb Q}$. To this datum one associates the quantum torus $T_q$, with $q=q(\theta)=e^{2i\pi\theta}$.
From the preceding subsection we see that $T_{q(\theta)}$ is simple if and only if $\theta\not\in{\mathbb Q}$, and that $T_{q(\theta)}$ and $T_{q(\theta')}$ are isomorphic if and only if  $\theta'=\theta+m$ or $\theta'=-\theta+m$ with $m\in{\mathbb Z}$.
This leads to the following.

\begin{sprop}
Consider two GWAs defined by  polynomials of degree two $a(h)=(h-\l_1)(h-\l_2)$ and $b(h)=(h-\mu_1)(h-\mu_2)$. Fix $\theta=\l_1-\l_2$, $\theta'=\mu_1-\mu_2$, $q=e^{2i\pi\theta}$, $q'=e^{2i\pi\theta'}$, and denote by ~$T_q$ and ~$T_{q'}$ the associated quantum tori. Then,
if ~$T_q$ (resp. $T_{q'}$) is simple then $A(a)$ (resp. $A(b)$) is simple and has finite global dimension. 

Assuming now that this condition holds for both $q$ and $q'$ in the following statements, then:
\begin{itemize}
 \item $A(a)\simeq A(b)$ if and only if $\theta=\pm\theta'$.
 \item $A(a)$ and $A(b)$ are Morita equivalent if and only if ~$T_q\simeq T_{q'}$.
\end{itemize}
\end{sprop}
\begin{demo}
 The first assertion and the first item are straightforward from previous remarks. For the last point, just note that $q'=q^{\pm 1}$ if and only if $\theta'=\pm\theta+m$, with $m\in{\mathbb Z}$.
\end{demo}

\begin{sremark}
{\rm The previous Proposition provides an alternative approach to Hodges' result concerning Morita equivalence for GWA when $n=2$.}

\end{sremark}

\medskip

The following subsection is devoted to obtain some generalisations in any degree.

\subsection{Rank $n$}

{\bf Notations.} For a polynomial $a(h)=\prod_{i=1}^n(h-\l_i)$ we will denote $\Theta(a)=(\theta_{ij})$ the matrix in $M_n({\mathbb Z})$ defined by $\theta_{ij}=\l_i-\l_j$. 
This matrix is not uniquely determined, since it actually depends on an indexing of the roots of $a$. In the sequel we will always assume that the polynomial $a$ is given with an indexing of its roots (counted with their multiplicities if necessary), and state our results up to a reindexing of these roots (see for example next Proposition).
Now we set $q_{ij}=e^{2i\pi\theta_{ij}}$ and  $Q(a)=(q_{ij})\in M_n({\mathbb C}^*)$. The matrix $Q(a)$ is multiplicatively antisymmetric, and $T_{Q(a)}$ will denote the quantum torus associated to these data.
We first note the following fact.

\begin{sprop}
 With the notations above, two generalized Weyl algebras $A(a)$ and $A(b)$ are isomorphic if and only if $\Theta(b)=\pm S^{-1}\Theta(a)S$ for a permutation matrix $S$.
\end{sprop}
\begin{demo}
 Denote by $(\l_1,\ldots,\l_n)$ and $(\mu_1,\ldots\mu_n)$ the roots of $a(h)$ and $b(h)$, counted with their multiplicities. It results from Theorem \ref{thmisogwa} that $A(a)$ and $A(b)$ are isomorphic if and only if there exists a permutation $\sigma$, a scalar $\beta$ and a sign $\epsilon$ such that $\mu_i=\epsilon\l_{\sigma i}+\beta$ for all $i$. From this one deduces easily the ``only if'' direction. For the reciprocal, assume that $\mu_i-\mu_j=\epsilon(\l_{\sigma i}-\l_{\sigma_j})$. Thanks to Theorem \ref{thmisogwa} one can assume that up to isomorphism $\epsilon=1$ and $\sigma={\rm Id}$, and that $\l_1=\mu_1=0$, and then $a(h)=b(h)$.
\end{demo}

\begin{scorol}
If the algebras $A(a)$ and $A(b)$ are isomorphic then the associated quantum tori $T_{Q(a)}$ and $T_{Q(b)}$ are isomorphic.
\end{scorol}
\begin{demo}
Permuting the generators with respect to the matrix $S$, one only has to prove that the matrix $Q(a)$ and its transpose define the same quantum torus.
According to the notations of Theorem \ref{thmqiso} the isomorphism is defined thanks to the matrix $M=(m_{ij})$
where $m_{11}=0$, $m_{1j}=1$ if $j\geq 2$; $m_{22}=0$, $m_{2j}=1$ if $j\neq 2$; $m_{ij}=-\delta_{ij}$ if $i\geq 2$. One uses the fact that $\l_{ij}\l_{jk}=\l_{ik}$ to verify that condition 2 of Theorem \ref{thmqiso} is satisfied. We leave the details to the reader.
\end{demo}


\begin{sremark}{\rm
This result strongly relies on the particular form of the parametrization  matrices we have here, and the fact that $\lambda_{ij}\lambda_{jk}=\lambda_{ik}$. For instance, taking $\l,\mu,\rho\in{\mathbb C}^*$ algebraically independant, the matrix
$\left({\small \begin{array}{ccc}
 1&\l&\mu\\ \l^{-1}&1&\rho\\ \mu^{-1}&\rho^{-1}&1
\end{array}}\right)$ and its transpose parametrize two quantum tori which are not isomorphic, since the corresponding matrix $G=(g_{ij})\in GL_3({\mathbb Z})$ in Theorem \ref{thmqiso} should satisfy $g_{11}^2g_{22}^2g_{33}^2=-1$.
}
\end{sremark}

We introduce now the following condition on $A(a)$.

\begin{sdefin}
 With the notations above, a generalized Weyl algebra $A(a)$ will be called \emph{$q$-simple} if the associated quantum torus $T_{Q(a)}$ is simple.
\end{sdefin}

\begin{sprop}
 Assume that the GWA $A(a)$ is $q$-simple. Then it is simple and has finite global dimension.
\end{sprop}
\begin{demo}
 By Proposition \ref{simple} we only have to show that if $\l_i-\l_j\in{\mathbb Z}$ then the matrix $\Theta(a)$ cannot satisfy condition (\ref{eqqsimple}). But this is clear by using the vector of ${\mathbb Z}^n$ with 1 in the $i$th place, $-1$ in the $j$th place and 
 $0$ everywhere else.
\end{demo}


\begin{sremark}
 {\rm In the case $n=2$, for a polynomial $a=(h-\l_1)(h-\l_2)$, being $q$-simple is equivalent to $\l_1-\l_2\not\in{\mathbb Q}$. This shows that  $q$-simplicity is strictly stronger  than simplicity and finite global dimension.}
\end{sremark}

Now we restate  condition 2. of Theorem \ref{thmqiso} in terms of matrices $\Theta(a)$ and $\Theta(b)$ associated to the roots of the polynomials $a$ and $b$. 

\begin{sprop}
 Let $a(h)=\prod_{i=1}^n(h-\l_i)$ and $b(h)=\prod_{i=1}^n(h-\mu_j)$ be two polynomials such that the GWAs $A(a)$ and $A(b)$ are $q$-simple. Then the quantum tori $T_{Q(a)}$ and $T_{Q(b)}$ are isomorphic if and only if there exist two matrices $M\in GL_n({\mathbb Z})$ and  $N\in M_n({\mathbb Z})$ such that $M^t\Theta(a)M=\Theta(b)+N$.
\end{sprop}
\qed

\medskip

It would be interesting  to relate this to Condition (\ref{cn}).
We end this discussion with some results in this direction concerning the case $n=3$.

\subsection{Case $n=3$}

The conditions above concerning the matrices can be restated, using cofactor matrices. More precisely,
let $\widehat{M}_{ij}$ be the matrix obtained from $M$ by deleting line $i$ and column $j$. Recall that if we denote by $cof(M)$ the matrix such that $cof(M)_{ij}=(-1)^{i+j}det(\widehat{M}_{ij})$ then $M\cdot cof(M)^t=det(M)\cdot Id$, so $det(cof(M))=1$ (since $n=3$ and $det(M)=\pm 1$).
We also have $det(cof(M))=det(cof(M))^t$, and $det(cof(M))=det(cof'(M))$, where $(cof'(M))_{ij}=det(\widehat{M}_{ij})$.
We rephrase in this case the conditions of the previous Proposition in terms of cofactor matrices. 

\begin{sprop}
Under the hypotheses of the above proposition, for $n=3$, the condition concerning matrices holds if and only if
\[cof'(M)^t\cdot \left(\begin{array}{c}
                      \lambda_{23}\\ \lambda_{13}\\ \lambda_{12}
                     \end{array}\right)
=\left(\begin{array}{c}
                      \mu_{23}\\ \mu_{13}\\ \mu_{12}
                     \end{array}\right)+\left(\begin{array}{c}
                      \gamma_{23}\\ \gamma_{13}\\ \gamma_{12}
                     \end{array}\right)\]
\end{sprop}

Taking into account that $\lambda_{12}=\lambda_{13}-\lambda_{23}$, and similarly for the $\mu$'s, we are able to 
establish a relation between Morita equivalences and isomorphisms of quantum tori for $n=3$.




\begin{stheoreme}
Fix $n=3$. If two generalized Weyl algebras $A(a)$ and $A(b)$ of degree $n$ are Morita equivalent, then their associated quantum tori are
isomorphic.
\end{stheoreme}

\begin{demo}
Given two Morita equivalent algebras $A(a)$ and $A(b)$, let $N=\left(\begin{array}{cc} n_1 &n_2\\ n_3&n_4
\end{array}\right)\in GL_2({\mathbb Z})$  be a matrix as in Section 3. We will construct a matrix  
$\widehat{N}\in GL_3({\mathbb Z})$  such that 
\[ \widehat{N}.\left(\begin{array}{c} \lambda_{23}\\ \lambda_{13}\\ \lambda_{12}\end{array}\right)=
\left(\begin{array}{c} \mu_{23}\\ \mu_{13}\\ \mu_{12}\end{array}\right) + \left(\begin{array}{c} \gamma_{23}\\ \gamma_{13}\\ \gamma_{12}\end{array}\right). \]

In fact, it is sufficient to take $\widehat{N}=\left(\begin{array}{ccc} n_4 &n_3 & 0 \\ n_2&n_1 & 0 \\ c & d& 1\end{array}\right)$,
where $c= n_1 -n_3 +1$ and $d= n_2 -n_4 -1$.

It is then straightforward to find a matrix $M\in GL_3({\mathbb Z})$ such that $\widehat{N}= cof'(M)^t$.

\end{demo}

\end{document}